\documentclass[conference]{IEEEtran}
\IEEEoverridecommandlockouts
\usepackage{cite}
\usepackage{amsmath,amssymb,amsfonts}
\usepackage{algorithmic}
\usepackage{graphicx}
\usepackage{textcomp}
\usepackage{xcolor}
\usepackage{dirtytalk}
\usepackage{comment}
\usepackage{float}
\usepackage{siunitx}
\def\BibTeX{{\rm B\kern-.05em{\sc i\kern-.025em b}\kern-.08em
    T\kern-.1667em\lower.7ex\hbox{E}\kern-.125emX}}

\begin{document}

\title{A data-driven method for quantifying the impact of a genetic circuit on its host  \\
\thanks{  \newline \textsuperscript{1} Department of Mechanical Engineering, University of California Santa Barbara, \textsuperscript{2} Pacific Northwest National Laboratories,  \textsuperscript{3} Department of Bioengineering, MIT,  \textsuperscript{4} Gingko Bioworks, \textsuperscript{5} Broad Foundry, MIT,  \textsuperscript{6} Texas Advanced Computing Center, University of Texas at Austin \newline \textsuperscript{*} To whom correspondence should be addressed: aqib@ucsb.edu, subhrajit.sinha@pnnl.gov, dorfany@mit.edu}
}
\author{\IEEEauthorblockN{Aqib Hasnain\textsuperscript{1*}, Subhrajit Sinha\textsuperscript{2*}, Yuval Dorfan\textsuperscript{3*}, Amin Espah Borujeni\textsuperscript{3}, Yongjin Park\textsuperscript{3}, Paul Maschhoff\textsuperscript{4},}
\IEEEauthorblockN{Uma Saxena\textsuperscript{5}, Joshua Urrutia\textsuperscript{6}, Niall Gaffney\textsuperscript{6}, Diveena Becker\textsuperscript{4}, Atsede Siba\textsuperscript{4}, Narendra Maheshri\textsuperscript{4},} 
\IEEEauthorblockN{Ben Gordon\textsuperscript{3,5}, Chris Voigt\textsuperscript{3}, and Enoch Yeung\textsuperscript{1}}
}

\maketitle

\begin{abstract}
Genetic circuits are designed to implement certain logic in living cells, keeping burden on the host cell minimal. However, manipulating the genome often will have a significant impact for various reasons (usage of the cell machinery to express new genes, toxicity of genes, interactions with native genes, etc.). In this work we utilize Koopman operator theory to construct data-driven models of transcriptomic-level dynamics from noisy and temporally sparse RNAseq measurements.  We show how Koopman models can be used to quantify impact on genetic circuits.  We consider an experimental example, using high-throughput RNAseq measurements collected from wild-type \textit{E. coli}, single gate components transformed in \textit{E. coli}, and a NAND circuit composed from individual gates in \textit{E. coli}, to explore how Koopman subspace functions encode increasing circuit interference on \textit{E. coli} chassis dynamics.  The algorithm provides a novel method for quantifying the impact of synthetic biological circuits on host-chassis dynamics.
\end{abstract}


\section{Introduction}

Synthetic biology is concerned with building useful biological circuits with predictable dynamics. Most genetic circuits depend on the local machinery of a host chassis. Since the host must expend energy in expressing exogenous genes, resources are taken from endogenous processes which can impact the function of both synthetic and native genes. Genetic circuit burden is often explicitly modeled as sequestration of cellular resources or competition for a limited set of binding sites.  These models are hypothesis driven, in that they rely on hypotheses or explicit knowledge of the nature of circuit burden.  

Functional genetic circuits require a minimal footprint on the host.  The higher the burden on the host, the more likely it is that the host mutates out the circuit. Even in cases where synthetic genes leave a minimal footprint on the host, an unpredictable change in dynamic behavior may occur due to the activation of synthetic genes. 
Rondelez showed that competition between synthetic and native genes have important effects on the global dynamics of the system \cite{rondelez2012competition}.  Recent studies involving host-circuit interactions looked at cross talk \cite{yeung2012quantifying} between genetic circuits and host resources for transcription \cite{gyorgy2015isocost}, translation \cite{borkowski2016overloaded,ceroni2015quantifying,gorochowski2016minimal,carbonell2015dealing,nystrom2018dynamic}, and protein degradation \cite{cookson2011queueing,qian2017resource}. These model-based approaches further our understanding of host-circuit interaction, however the models are often based on biophysical mechanisms that are difficult to validate or observe. Transcriptomics and proteomics resolve the activity of thousands of genes, providing a rich resource for learning models to answer key questions without the need for hypothesis-driven modeling. How can these measurements be leveraged through a data-driven approach to better understand host-circuit interaction and genetic stability? 

Spectral methods have been increasingly popular in the data-driven analysis of nonlinear dynamical systems. Recently, researchers working in Koopman operator theory
have shown that it is possible to identify and learn the fundamental modes of a nonlinear dynamical system from data \cite{rowley_mezic_bagheri_schlatter_henningson_2009}. The seminal work by Schmid in developing dynamic mode decomposition (DMD) has led to an enormous growth in the use of Koopman spectral analysis of nonlinear dynamical systems \cite{schmid2010dynamic}. More recently, learning higher dimensional Koopman operators from data has become computationally tractable, largely due to advances in integrating machine learning and deep learning to generate efficient representations of observable bases \cite{yeung2017learning,lusch2018deep,otto2019linearly}. Often in biology and especially in omics measurements, the data are temporally sparse. Sinha and Yeung developed a method for computing the Koopman operator from sparse data \cite{sinha_yeung_2019}. 


Synthetic biological circuit design is often viewed from a reductionist's perspective. Biological parts are designed and optimized for modularity, so that composition gives rise to predictable behavior. The challenge is that composition, while at times successfully gives rise to predictable {\it observed} behavior, has an unknown {\it emergent} effect on the host, and by the principles of feedback, the circuit as well \cite{cardinale2012contextualizing}.  

In this paper we develop a completely novel algorithm, \textit{structured DMD}, to complement bottom-up genetic circuit design approaches in synthetic biology. Structured DMD is a purely data-driven model discovery framework that takes advantage of a part-by-part construction process to decouple emergent phenomena from isolated part dynamics.  It reduces the total model complexity in the model identification process by adopting a hierarchical approach to identifying components of the model in stages.   The decomposition we obtain is additive, due to nice linear mathematical properties endowed by Koopman operators \cite{mezic2005spectral}.  We showcase our algorithm on a coupled oscillator system, but then consider a real genetic circuit design problem using a NAND gate designed from TetR orthologs \cite{stanton2014genomic}.  Full-state but temporally sparse RNAseq measurements collected from wild-type \textit{E. coli}, single gate components transformed in \textit{E. coli}, and a NAND circuit composed from individual gates in \textit{E. coli} are used to explore how Koopman subspace functions encode increasing circuit interference on \textit{E. coli} chassis dynamics.

\section{Koopman operator formulation and Dynamic mode decomposition} \label{sec:Koop}
We briefly introduce Koopman operator theory (see \cite{mezic2005spectral} for a full discussion); as we will use it throughout this paper. Consider a discrete time open-loop nonlinear system of the form 
\begin{equation}
 x_{t+1} =f(x_t)\\
 \label{eq:sys}
\end{equation}
where $f: M \subset \mathbb{R}^n \rightarrow M$ is an analytic vector field on the state space. The Koopman operator of (\ref{eq:sys}), $\mathcal{K}$ : $\mathcal{F}$ $\rightarrow$ $\mathcal{F}$, is a linear operator that acts on observable functions $\psi (x_k)$ and propagates them forward in time as 
\begin{equation}
    \psi (x_{t+1})=\mathcal{K}\psi (x_t).
    \label{eq:KoopEqInf}
\end{equation}
Here ${\mathcal F}$ is the space of observable functions that is invariant under the action of $\mathcal{K}$.

Using data-driven approaches, commonly DMD  \cite{schmid2010dynamic} or extended DMD \cite{williams_kevrekidis_rowley_2015}, an approximation to the Koopman operator, $K$, can be computed. The approach taken to compute an approximation to the Koopman operator in both DMD and extended DMD is to solve the following optimization problem 
\begin{equation}
\min_{K} || \Psi(X_f) - K\Psi(X_p)||
\label{eq:learnKoop}
\end{equation}
where 
$  X_f \equiv \begin{bmatrix} x_{1} & \hdots & x_{N-1} \end{bmatrix},$ $ X_p  \equiv \begin{bmatrix} x_{2} & \hdots & x_{N} \end{bmatrix}
$
are snapshot matrices formed from the discrete-time dynamical system (\ref{eq:sys}) and 
$
  \Psi(X) \equiv \begin{bmatrix} \psi_1(x) & \hdots & \psi_R(x) \end{bmatrix}
$
is the mapping from physical space into the space of observables. DMD is a special case of extended DMD where $\psi(x) = x$. It was shown by Rowley et al. that the approximate Koopman operator obtained from DMD is closely related to a spectral analysis of the linear but infinite-dimensional Koopman operator \cite{rowley_mezic_bagheri_schlatter_henningson_2009}.

\section{Structured dynamic mode decomposition} \label{sec:sdmd}

In this section we, for the first time, introduce structured dynamic mode decomposition (structured DMD). The structured DMD algorithm takes advantage of bottom-up design approaches where an original system is built upon by adding parts layer by layer to achieve complex dynamical behaviors. 

The Koopman model (\ref{eq:KoopEqInf}) can be decomposed into original and added equations, or its parts in the bottom-up design approach, written as
\begin{equation}
    \begin{bmatrix} \psi_O(x_{t+1})\\\psi_A(x_{t+1}) \end{bmatrix} \equiv \begin{bmatrix} K_{OO} & K_{OA} \\ K_{AO} & K_{AA} \end{bmatrix} \begin{bmatrix} \psi_O(x_{t})\\\psi_A(x_{t}) \end{bmatrix}
    \label{eq:KoopEqOA}
\end{equation}
where the subscripts $O$ and $A$ correspond to the original and added components. The matrix of Koopman operators in (\ref{eq:KoopEqOA}) is unknown, but can be obtained using the standard techniques outlined in Section \ref{sec:Koop}. However, these approaches would not allow the decoupling of the underlying dynamics of the original components from added  components. Therefore, it would not be possible to determine the impact that new components have on the original components. With this algorithm, we propose to discover the underlying (original) dynamics directly from data first, and then subsequently learn the interaction dynamics as an additive perturbation in the Koopman model.

If we want to solely understand the impact of added components on the original components of a system, we first learn the original-original interaction dynamics matrix $K_{OO}$ from the original system without added parts. This original system has Koopman model
\begin{equation}
    \psi_O(x_{t+1}) = K_{OO}\psi_O(x_t)
\end{equation}
where $K_{OO}$ is learned through the optimization problem (\ref{eq:learnKoop}). The original-added interaction dynamics matrix $K_{OA}$ can now be learned by viewing it as an added perturbation in the original Koopman model, i.e. the new model is
\begin{equation}
    \psi_O(x_{t+1}) = \underbrace{K_{OO}}_{known}\psi_O(x_t) + K_{OA}\psi_A(x_t). 
\end{equation}
Here we have already learned $K_{OO}$ and of interest to us is $K_{OA}$ which can now be learned directly from data. In this way, we can completely decouple the underlying original dynamics from the effect that any added parts have on the original dynamics. 

\section{Impact of genetic circuit on host} \label{sec:impact}
\begin{figure*}
    \centering
    \vspace{-10mm}
    \includegraphics[width=1.9\columnwidth]{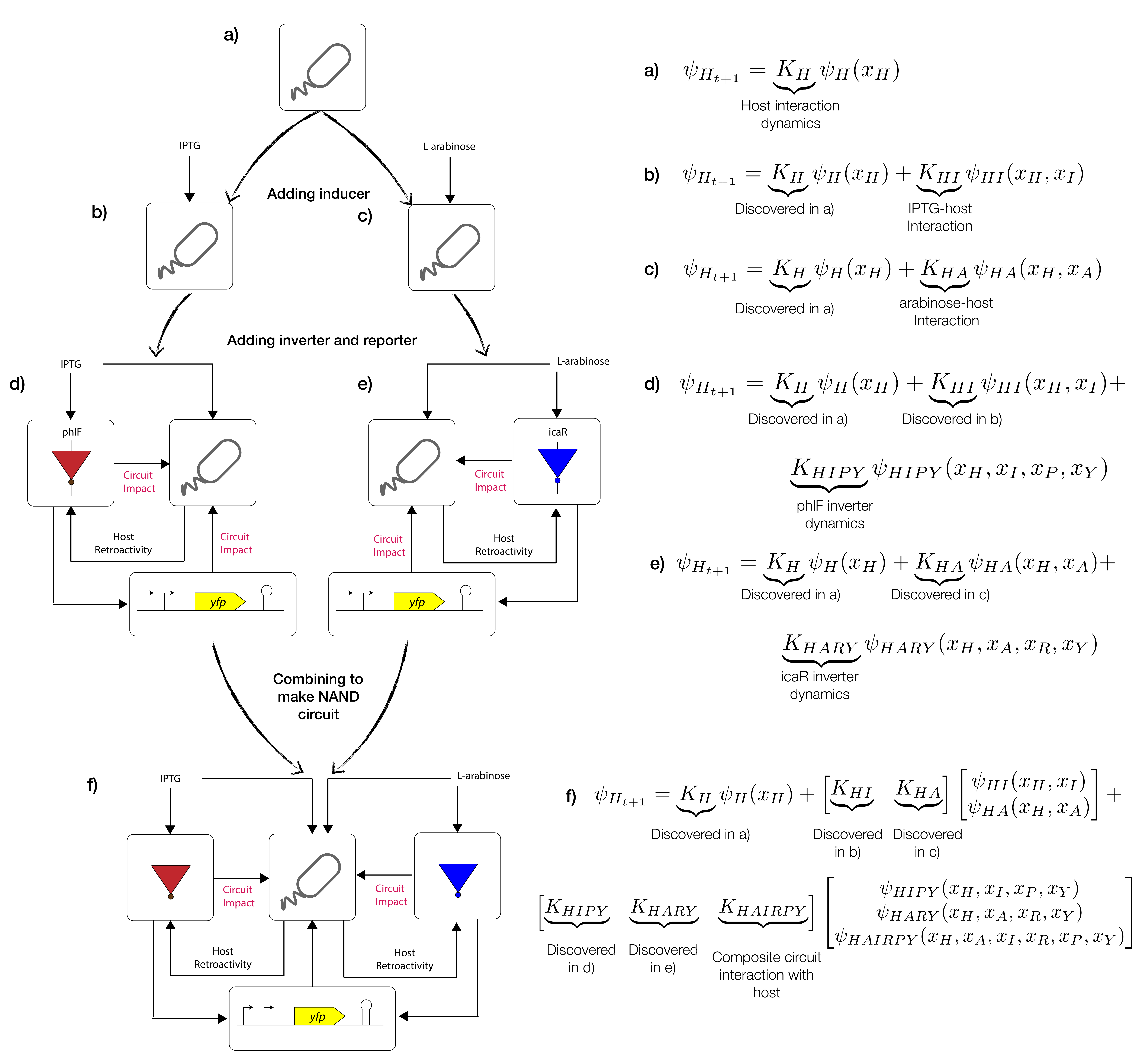}
    \caption{Schematic of the bottom-up design of a NAND gate in \textit{E. coli}. a) \textit{E. coli}, b) \textit{E. coli} with IPTG input c) \textit{E. coli} with L-arabinose input d) \textit{E. coli} with IPTG input, PhlF gate, and YFP reporter e) \textit{E. coli} with L-arabinose input, IcaR gate, and YFP reporter, f) complete NAND circuit. Under each design iteration is the associated host Koopman model.}
    
    \label{fig:CircuitParts}
\end{figure*}
First, note that typically RNAseq measurements are sparse in time (two timepoints and four replicates in this case). Sinha and Yeung \cite{sinha_yeung_2019} have addressed the problem of computation of the Koopman operator when the data is sparse. For each data tuple $(x_i,x_{i+1})$, the artificial data point $(x_i + \delta x_i, x_{i+1} + \delta x_{i+1})$ is added. Artificial snapshot matrices are formed as 
\begin{equation*}
\begin{aligned}
  X_f &\equiv \begin{bmatrix} x_{1} & \hdots & x_{N-1} & x_{1}+\delta x_{1} & \hdots & x_{N-1}+\delta x_{N-1} \end{bmatrix},  \\ X_p & \equiv \begin{bmatrix} x_{2} & \hdots & x_{N} & x_2+\delta x_2 & \hdots & x_{N} +\delta x_{N} \end{bmatrix}
\end{aligned}
\end{equation*}
These artificial data points (which are sufficiently small perturbations) are added to the sparse data set to enrich the data. Robust optimization-based techniques are then used to compute the approximate Koopman operator. The optimization problem to be solved is 
\begin{equation*}
    \min_{K} || \Psi(X_f) - K\Psi(X_p) ||_F + \lambda || K ||_F
\end{equation*}
where $\lambda$ is a regularization parameter. 

\begin{figure}
    \centering
    \includegraphics[width=.7\columnwidth]{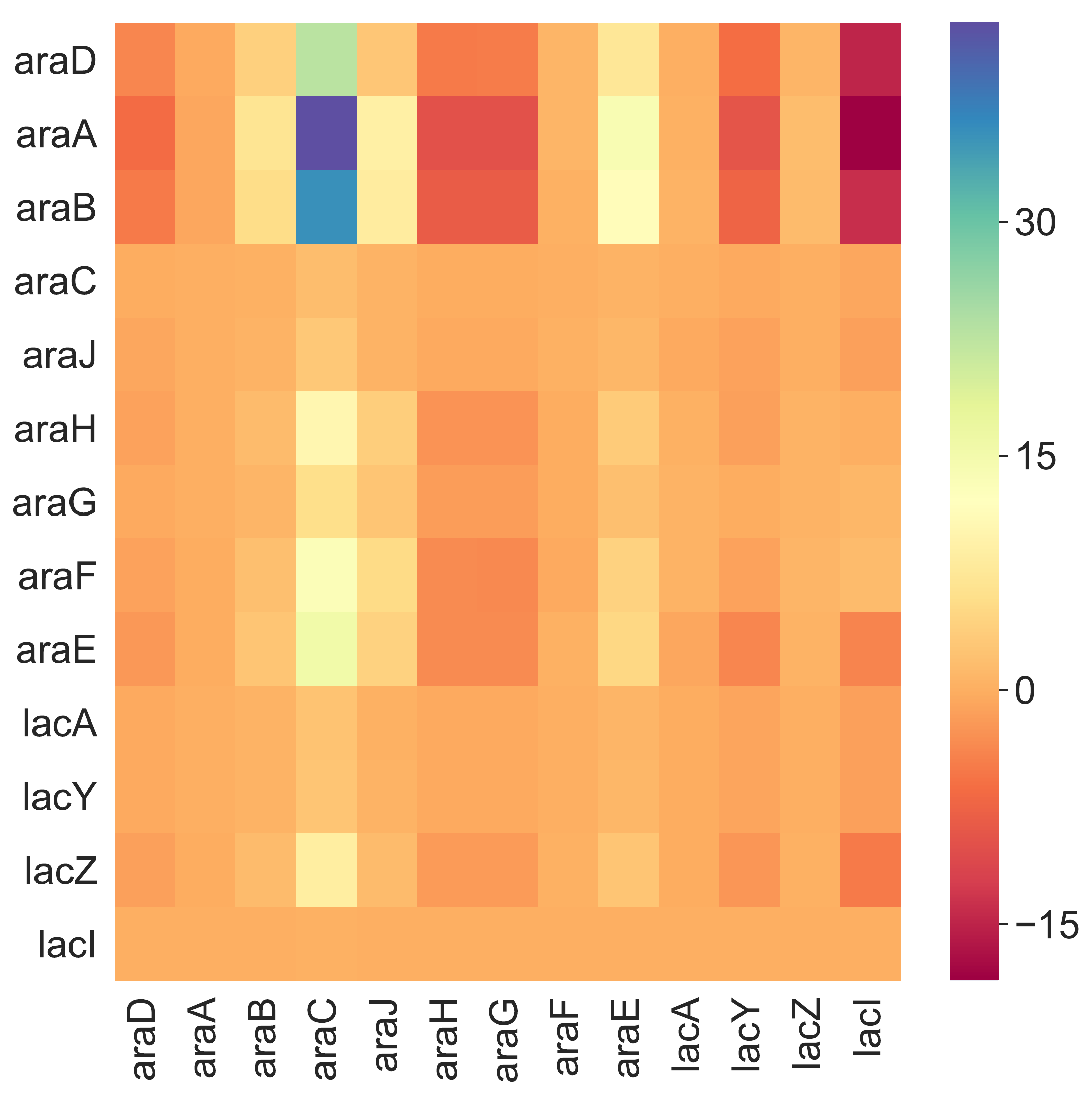}
    \caption{The Koopman operator estimated from structured DMD of the host dynamics in response to arabinose and IPTG induction. The color scale represents likely causal interaction; positive causal interaction represented by positive values and negative causal interaction represented by negative values.}
    \label{fig:Kol}
    \vspace{-4mm}
\end{figure}

\begin{figure}
    \centering
    \includegraphics[width=\columnwidth]{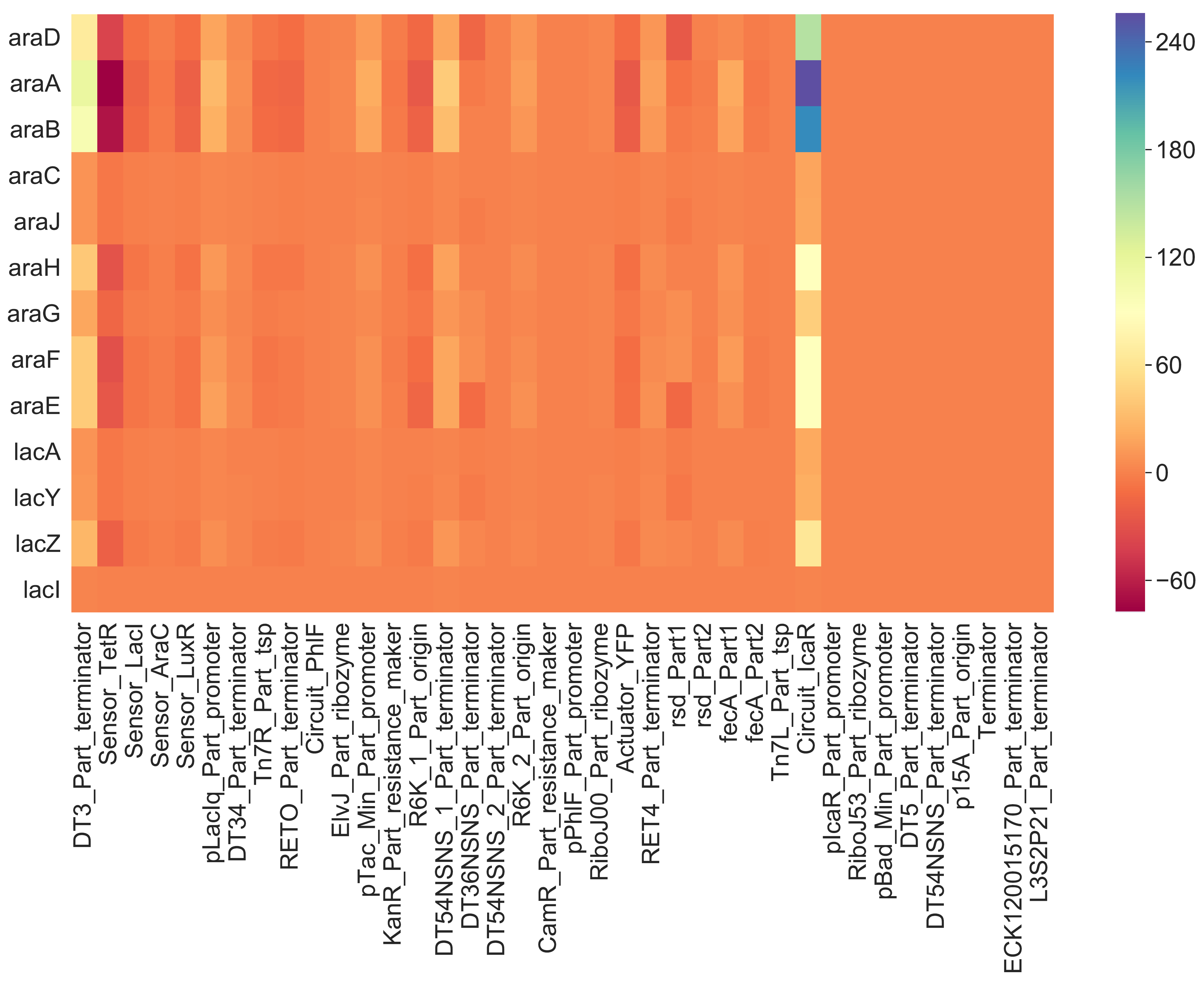}
    \vspace{-4mm}
    \caption{The input Koopman operator estimated from structured DMD modeling circuit-to-host interaction of the NAND circuit. The color scale represents likely causal interaction; positive causal interaction represented by positive values and negative causal interaction represented by negative values.}
    \label{fig:K_NAND}
\end{figure}

Figure \ref{fig:CircuitParts} shows the bottom-up construction of a NAND gate in \textit{E. coli} where each iteration also has an associated host Koopman model (under each schematic). The composite interaction matrix $K_{HAIRPY}$ defines the impact of the NAND gate on the host. To learn this matrix, we first learn the underlying host dynamics $K_H$ in figure \ref{fig:CircuitParts}a. A heatmap of $K_H$ can be seen in figure \ref{fig:Kol}. Inducer-host interactions $K_{HI}$ and $K_{HA}$ in figure \ref{fig:CircuitParts}b and \ref{fig:CircuitParts}c are computed next. The subscripts $I$ and $A$ correspond to the inducers IPTG and L-arabinose, respectively. PhlF and IcaR inverters are then added along with \textit{yfp} reporters as seen in figure \ref{fig:CircuitParts}d and \ref{fig:CircuitParts}e. $K_{HIPY}$ and $K_{HARY}$ are learned from these systems where the subscripts $P$, $R$, and $Y$ denote the PhlF inverter, IcaR inverter, and \textit{yfp} reporters, respectively. At this stage, we can learn $K_{HAIRPY}$ since all the terms in the Koopman model of figure \ref{fig:CircuitParts}f are now known. Figure \ref{fig:K_NAND} shows a heatmap of $K_{HAIRPS}$. 

The discovered host Koopman model from wild type MG1655K12 {\it E. coli} reveals a fundamentally antagonistic relationship between the LacI and AraC control modules in the transition from the first time point (log phase) to the second timepoint (stationary phase).   We see that arabinose induction activates most of the Ara operon genes, as expected, but simultaneously creates a negative inhibitory effect on LacY, LacZ, and LacI.  Conversely, induction with IPTG has a significant downregulating effect on the Ara operon genes, specifically the cluster of genes downstream of the pBAD cluster.  This is consistent with prior analysis of the hierarchy of sugar utilization \cite{aidelberg2014hierarchy}.  Crosstalk mechanisms mediated by the catabolite repression protein (CRP) and cyclic AMP pathway prioritize Lac operon activity over Ara operon activity when both sugars are present.  While IPTG is not a sugar, it acts as a structural analog and thus induces the same diauxic response.   We have verified this phenomena, using a data-driven approach, using only two timepoints and four biological replicates from noisy RNAseq measurements. 

The import of this finding is that the NAND circuit inherently activates the diauxic response mechanism to its advantage.  Since PhlF and IcaR are designed and intended to act independently, mutual repression of their underlying host machinery results in stronger underlying XOR logic and mutual coupling.  That is, activation with arabinose will result in indirect repression of the lactose operon, and vice-versa.   

Finally, we discovered that IcaR gene expression induces a positive-feedback loop with the AraBAD cluster. This in turn, results in elevated IcaR expression, which induces cytotoxicity.   When calculating the Frobenius norm, as a total sum measure of circuit-to-host impact, we found that arabinose induction in the host had an impact of $ \frac{||K_{HA}||_F}{p\times q} = \num{4e-7} $ ($p$ and $q$ are the dimensions of $K_{HA}$), while induction of the IcaR component had an impact of $ \frac{||K_{HIPY}||_F}{r\times s} = \num{2e-2}$, nearly 5 orders of magnitude greater.  Even though PhlF and IcaR had comparable per-term gain over the 429 genes we analyzed, IcaR impacted 420 genes, while PhlF only impacted 269 genes.    When analyzing DNA sequencing data, we found that the IcaR gene had been deleted from the NAND circuit on the genome; the disparity in gain and widespread influence between the IcaR and the other circuit components provides a hypothesis for IcaR mutation.  The IcaR part imposes considerable widespread perturbation on host genes, which is evidence of cytotoxicity leading to mutation.    

 
\section*{Acknowledgements}
The authors gratefully acknowledge the funding of DARPA grants FA8750-17-C-0229, HR001117C0092, HR001117C0094, DEAC0576RL01830. The authors would also like to thank Professors Igor Mezic, Alexandre Mauroy, Nathan Kutz, Steve Haase, John Harer, and Eric Klavins for insightful discussions.  
Any opinions, findings, conclusions, or recommendations expressed in this material are those of the authors and do not necessarily reflect the views of the Defense Advanced Research Project Agency, the Department of Defense, or the United States government. This material is based on work supported by DARPA and AFRL under contract numbers FA8750-17-C-0229, HR001117C0092, HR001117C0094, DEAC0576RL01830.
\bibliographystyle{abbrv}
\bibliography{main}

\end{document}